\documentclass[preprint,12pt,authoryear]{elsarticle}
\usepackage[]{graphicx,float,latexsym,times}
\usepackage{amsfonts,amstext,amsmath,amssymb,amsthm}
\usepackage{caption2}

\newcommand{\Var}{\mbox{Var}}

\newcommand{\be}{\mbox{\boldmath {$e$}}}

\newcommand{\bh}{\mbox{\boldmath {$h$}}}
\newcommand{\bH}{\mbox{\boldmath {$H$}}}

\newcommand{\bI}{\mbox{\boldmath {$I$}}}

\newcommand{\bO}{\mbox{\boldmath {$O$}}}

\newcommand{\bP}{\mbox{\boldmath {$P$}}}

\newcommand{\bS}{\mbox{\boldmath {$S$}}}

\newcommand{\bu}{\mbox{\boldmath {$u$}}}

\newcommand{\bx}{\mbox{\boldmath {$x$}}}
\newcommand{\bX}{\mbox{\boldmath {$X$}}}

\newcommand{\bz}{\mbox{\boldmath {$z$}}}
\newcommand{\bZ}{\mbox{\boldmath {$Z$}}}
\newcommand{\bze}{\mbox{\boldmath {$0$}}}
\newcommand{\bone}{\mbox{\boldmath {$1$}}}

\newcommand{\bmu}{\mbox{\boldmath $ \mu $}}
\newcommand{\bSig}{\mbox{\boldmath $ \Sigma $}}

\newcommand{\bSigma}{\mbox{\boldmath $ \Sigma $}}

\newcommand{\bLam}{\mbox{\boldmath $ \Lambda $}}

\newcommand{\tr}{\mbox{tr}}
\newcommand{\argmin}{\mathop{\rm argmin}\limits}

\theoremstyle{plain} 
\newtheorem{theorem}{Theorem}[section]
\newtheorem{lemma}{Lemma}[section]
\newtheorem{corollary}{Corollary}[section]
\newtheorem{proposition}{Proposition}[section]

\theoremstyle{definition} 

\newtheorem{remark}{Remark}[section]

\numberwithin{equation}{section} 

\begin{document}

\begin{frontmatter}



\title{Asymptotic properties of the first principal component
and equality tests of covariance matrices 
in high-dimension, low-sample-size context}


\author[a]{Aki Ishii}
\author[b]{Kazuyoshi Yata}
\author[b]{Makoto Aoshima\fnref{fn1}}
\fntext[fn1]{Institute of Mathematics, University of Tsukuba, Ibaraki 305-8571, Japan;\\ 
\hspace{6mm}Fax: +81-298-53-6501}
\ead{aoshima@math.tsukuba.ac.jp}

\address[a]{Graduate School of Pure and Applied Sciences, University of Tsukuba, Ibaraki, Japan}
\address[b]{Institute of Mathematics, University of Tsukuba, Ibaraki, Japan}

\begin{abstract}
A common feature of high-dimensional data is that the data dimension is high, however, the sample size is relatively low. 
We call such data HDLSS data.
In this paper, we study asymptotic properties of the first principal component in the HDLSS context and apply them to equality tests of covariance matrices for high-dimensional data sets.
We consider HDLSS asymptotic theories as the dimension grows for both the cases when the sample size is fixed and the sample size goes to infinity.
We introduce an eigenvalue estimator by the noise-reduction methodology and provide asymptotic distributions of the largest eigenvalue in the HDLSS context.
We construct a confidence interval of the first contribution ratio. 
We give asymptotic properties both for the first PC direction and PC score as well. 
We apply the findings to equality tests of two covariance matrices in the HDLSS context.  
We provide numerical results and discussions about the performances both on the estimates of the first PC and the equality tests of two covariance matrices. 
\end{abstract}


\end{frontmatter}


\section{Introduction}
\label{s:Intro}
One of the features of modern data is the data dimension $d$ is high and the sample size $n$ is relatively low.
We call such data HDLSS data.
In HDLSS situations such as $d/n\to\infty$, new theories and methodologies are required to develop for statistical inference based on the large sample theory.
One of the approaches is to study geometric representations of HDLSS data and investigate the possibilities to make use of them in HDLSS statistical inference.
Hall et al. (2005), Ahn et al. (2007), and Yata and Aoshima (2012) found several conspicuous geometric descriptions of HDLSS data when $d\to\infty$ while $n$ is fixed. 
The HDLSS asymptotic studies usually assume either the normality as the population distribution or a $\rho$-mixing condition as the dependency of random variables in a sphered data matrix.
See Jung and Marron (2009) and Jung et al. (2012).
However, Yata and Aoshima (2009) developed an HDLSS asymptotic theory without assuming those assumptions and showed that the conventional principal component analysis (PCA) cannot give consistent estimation in the HDLSS context.   
In order to overcome this inconvenience, Yata and Aoshima (2012) provided the {\it noise-reduction (NR) methodology} that can successfully give consistent estimators of both the eigenvalues and eigenvectors together with the principal component (PC) scores.
Furthermore, Yata and Aoshima (2010, 2013) created the {\it cross-data-matrix (CDM) methodology} that is a nonparametric method to ensure consistent estimation of those quantities.    
Given this background, Aoshima and Yata (2011, 2013) developed a variety of inference for HDLSS data such as given-bandwidth confidence region, two-sample test, test of equality of two covariance matrices, classification, variable selection, regression, pathway analysis and so on along with the sample size determination to ensure prespecified accuracy for each inference.

In this paper, suppose we have a $d\times n$ data matrix, $\bX_{(d)}=[\bx_{1(d)},...,\bx_{n(d)}]$, where $\bx_{j(d)}=(x_{1j(d)},...,$
$x_{dj(d)})^T,\ j= 1,...,n$, are independent and identically distributed (i.i.d.) as a $d$-dimensional distribution with a mean vector $\bmu_d$ and covariance matrix $\bSig_d\ (\ge \bO)$.
We assume $n\ge 3$. 
The eigen-decomposition of $\bSig_d$ is given by $\bSig_d=\bH_d\bLam_d\bH_d^T$, where $\bLam_d$ is a diagonal matrix of eigenvalues, $\lambda_{1(d)}\ge \cdots \ge \lambda_{d(d)}(\ge 0)$, and $\bH_d=[\bh_{1(d)},...,\bh_{d(d)}]$ is an orthogonal matrix of the corresponding eigenvectors. 
Let $ \bX_{(d)}-[\bmu_d,...,\bmu_d]=\bH_d\bLam_d^{1/2}\bZ_{(d)}$. 
Then, $\bZ_{(d)}$ is a $d\times n$ sphered data matrix from a distribution with the zero mean and the identity covariance matrix. 
Here, we write $\bZ_{(d)}=[\bz_{1(d)},...,\bz_{d(d)}]^T$ and $\bz_{j(d)}=(z_{j1(d)},...,z_{jn(d)})^T,\ j=1,...,d$. 
Note that $E(z_{ji(d)}z_{j'i(d)})=0\ (j\neq j')$ and $\Var(\bz_{j(d)})=\bI_n$, where $\bI_n$ is the $n$-dimensional identity matrix. 
The $i$-th true PC score of $\bx_{j(d)}$ is given by $\bh_{i(d)}^T(\bx_{j(d)}-\bmu_d)=\lambda_{i(d)}^{1/2} z_{ij(d)}$ (hereafter called $s_{ij(d)}$). 
Note that Var$(s_{ij(d)})=\lambda_{i(d)}$ for all $i,j$. 
Hereafter, the subscript $d$ will be omitted for the sake of simplicity when it does not cause any confusion. 
We assume that $\lambda_1$ has multiplicity one in the sense that $\liminf_{d\to \infty}\lambda_{1}/\lambda_{2}>1$. 
Also, we assume that $\limsup_{d\to \infty} E(z_{ij}^4)<\infty$ for all $i,j$ and $P( \lim_{d\to \infty} ||\bz_1|| \neq 0 )=1$. 
Note that if $\bX$ is Gaussian, $z_{ij}$s are i.i.d. as the standard normal distribution, $N(0,1)$. 
As necessary, we consider the following assumption for the normalized first PC scores, $z_{1j}\ (=s_{1j}/\lambda_1^{1/2})$, $j=1,...,n$: 
\begin{description}
  \item[ (A-i)]\quad $z_{1j},\ j=1,...,n,$ are i.i.d. as $N(0,1)$.
\end{description}
Note that $P( \lim_{d\to \infty} ||\bz_1|| \neq 0 )=1$ under (A-i). 
Let us write the sample covariance matrix as $\bS=(n-1)^{-1}(\bX-\overline{\bX})(\bX-\overline{\bX})^T =(n-1)^{-1}\sum_{j=1}^n(\bx_j-\bar{\bx})(\bx_j-\bar{\bx})^T$, where $\overline{\bX}=[\bar{\bx},...,\bar{\bx}]$ and $\bar{\bx}=\sum_{j=1}^n\bx_j/n$. 
Then, we define the $n \times n$ dual sample covariance matrix by $\bS_{D}=(n-1)^{-1}(\bX-\overline{\bX})^T(\bX-\overline{\bX})$.
Let $\hat{\lambda}_{1}\ge\cdots\ge\hat{\lambda}_{n-1}\ge 0$ be the eigenvalues of $\bS_{D}$.  
Let us write the eigen-decomposition of $\bS_{D}$ as $\bS_{D}=\sum_{j=1}^{n-1}\hat{\lambda}_{j}\hat{\bu}_{j}\hat{\bu}_{j}^T $, where $\hat{\bu}_j=(\hat{u}_{j1},...,\hat{u}_{jn})^T$ denotes a unit eigenvector corresponding to $\hat{\lambda}_{j}$. 
Note that $\bS$ and $\bS_D$ share non-zero eigenvalues.

In this paper, we study asymptotic properties of the first principal component in the HDLSS context and apply them to equality tests of covariance matrices for high-dimensional data sets.
We consider HDLSS asymptotic theories as $d\to \infty$ for both the cases when $n$ is fixed and $n\to \infty$.
In Section 2, we introduce an eigenvalue estimator by the NR methodology and provide asymptotic distributions of the largest eigenvalue in the HDLSS context.
We construct a confidence interval of the first contribution ratio. 
In Section 3, we give asymptotic properties both for the first PC direction and PC score as well. 
In Section 4, we apply the findings to equality tests of two covariance matrices in the HDLSS context.  
Finally, in Section 5, we provide numerical results and discussions about the performances both on the estimates of the first PC and the equality tests of two covariance matrices. 
\section{Largest eigenvalue and its contribution rate}
In this section, we give asymptotic distributions of the largest eigenvalue and construct a confidence interval of the first contribution rate. 
\subsection{Asymptotic distributions of the largest eigenvalue}
Let $\delta_i=\tr(\bSig^2)-\sum_{s=1}^i\lambda_s^2=\sum_{s=i+1}^d\lambda_s^2$ for $i=1,...,d-1$. 
We consider the following assumptions for the largest eigenvalue: 
\begin{description}
  \item[ (A-ii)]\quad  $\displaystyle \frac{\delta_1}{\lambda_1^2}=o(1)$ as $d\to \infty$ when $n$ is fixed;
$\displaystyle \frac{\delta_{i_*}}{\lambda_1^2}=o(1)$ as $d\to \infty$ for some fixed $i_*\ (<d)$ when $n\to \infty$. 
  \item[ (A-iii)]\quad $\displaystyle \frac{\sum_{r,s\ge 2}^d \lambda_r\lambda_{s}E\{(z_{rk}^2-1)(z_{sk}^2-1)\}}{n\lambda_1^2}=o(1)$ as $d\to \infty$ either when $n$ is fixed or $n\to \infty$.
\end{description}
Note that (A-iii) holds when $\bX$ is Gaussian and (A-ii) is met. 
Let ${\bz}_{oj}=\bz_{j}-(\bar{z}_{j},...,\bar{z}_{j})^T,\ j=1,...,p$, where $\bar{z}_{j}=n_{}^{-1}\sum_{k=1}^{n}z_{jk}$. 
Let $\kappa=\tr(\bSig)-\lambda_1=\sum_{s=2}^d\lambda_s$. 
Then, we have the following result.
\begin{proposition}
Under (A-ii) and (A-iii), it holds that 
$$
\frac{\hat{\lambda}_1}{\lambda_1}-||{\bz}_{o1}/\sqrt{n-1}||^2- \frac{ \kappa}{\lambda_1(n-1)}=o_p(1)
$$
as $d\to \infty$ either when $n$ is fixed or $n\to \infty$.
\end{proposition}
\begin{remark}
Jung et al. (2012) gave a result similar to Proposition 2.1 when $\bX$ is Gaussian, $\bmu=\bze$ and $n$ is fixed.
\end{remark}

It holds that $E(||{\bz}_{o1}/\sqrt{n-1}||^2)=1$ and $||{\bz}_{o1}/\sqrt{n-1}||^2=1+o_p(1)$ as $n\to \infty$.
If $\kappa/(n\lambda_1)=o(1)$ as $d\to\infty$ and $n\to\infty$, $\hat{\lambda}_1$ is a consistent estimator of $\lambda_1$.  
When $n$ is fixed, the condition `$\kappa/\lambda_1=o(1)$' is equivalent to `$\lambda_1/\tr(\bSig)=1+o(1)$' in which the contribution ratio of the first principal component is asymptotically $1$. 
In that sense, `$\kappa/\lambda_1=o(1)$' is quite strict condition in real high-dimensional data analyses. 
Hereafter, we assume $\liminf_{d\to \infty} \kappa/\lambda_1>0$. 

Yata and Aoshima (2012) proposed a method for eigenvalue estimation called the {\it noise-reduction (NR) methodology} that was brought by a geometric representation of $\bS_D$.
If one applies the NR methodology to the present case, $\lambda_i$s are estimated by
\begin{equation}
\tilde{\lambda}_{i}=\hat{\lambda}_{i}-\frac{\tr(\bS_{D})-\sum_{j=1}^i\hat{\lambda}_{j} }{n-1-i}\quad (i=1,...,n-2).
\label{2.1}
\end{equation}
Note that $\tilde{\lambda}_i\ge 0$ w.p.1 for $i=1,...,n-2$. 
Also, note that the second term in (\ref{2.1}) with $i=1$ is an estimator of $\kappa/(n-1)$. 
See Lemma 2.1 in Section 2.2 for the details.
Yata and Aoshima (2012, 2013) showed that $\tilde{\lambda}_i$ has several consistency properties when $d\to \infty$ and $n\to \infty$. 
On the other hand, Ishii et al. (2014) gave asymptotic properties of $\tilde{\lambda}_1$ when $d\to \infty$ while $n$ is fixed. 
The following theorem summarizes their findings: 
\begin{theorem}
Under (A-ii) and (A-iii), it holds that as $d\to \infty$ 
\begin{align*}
&\frac{\tilde{\lambda}_1}{\lambda_1}= \left\{ \begin{array}{ll}
||{\bz}_{o1}/\sqrt{n-1}||^2+o_p(1) & \mbox{ when $n$ is fixed}, \\
& \\[-3mm]
1+o_p(1) & \mbox{ when $n\to \infty$}. \\
\end{array} \right.
\end{align*}
Under (A-i) to (A-iii), it holds that as $d\to \infty$ 
\begin{align*}
&(n-1)\frac{\tilde{\lambda}_1}{\lambda_1}\Rightarrow \chi_{n-1}^2 \hspace{2.3cm} \mbox{when $n$ is fixed},\\
&\sqrt{\frac{n-1}{2}}\Big(\frac{\tilde{\lambda}_1}{\lambda_1}-1 \Big) \Rightarrow N(0,1)\quad \mbox{when $n\to \infty$.}
\end{align*}
Here, $``\Rightarrow"$ denotes the convergence in distribution and $\chi_{n-1}^2$ denotes a random variable distributed as $\chi^2$ distribution with $n-1$ degrees of freedom.  
\end{theorem}
\subsection{Confidence interval of the first contribution ratio}
We consider a confidence interval for the contribution ratio of the first principal component. 
Let $a$ and $b$ be constants satisfying $P(a \le \chi_{n-1}^2 \le b )=1-\alpha$, where $\alpha\in (0,1)$.
Then, from Theorem 2.1, under (A-i) to (A-iii), it holds that
\begin{align}
&P\Big(\frac{\lambda_1}{\tr(\bSigma)}\in \Big[\frac{(n-1)\tilde{\lambda}_1}{b\kappa+(n-1)\tilde{\lambda}_1},\frac{(n-1)\tilde{\lambda}_1}{a\kappa+(n-1)\tilde{\lambda}_1}\Big]\Big)\notag \\
&=P\Big(a\le (n-1) \frac{\tilde{\lambda}_1}{\lambda_1} \le b \Big)=1-\alpha+o(1) 
\label{2.2}
\end{align}
as $d\to \infty$ when $n$ is fixed.
We need to estimate $\kappa$ in (\ref{2.2}). 
Here, we give a consistent estimator of $\kappa$ by $\tilde{\kappa}=(n-1)(\tr(\bS_D)-\hat{\lambda}_1)/(n-2)=\tr(\bS_D)-\tilde{\lambda}_1$. 
Then, we have the following results.
\begin{lemma}
Under (A-ii) and (A-iii), it holds that 
$$
\frac{\tilde{\kappa}}{\kappa}=1+o_p(1)\quad \mbox{and}\quad \frac{\tilde{\kappa}}{\lambda_1}=\frac{\kappa}{\lambda_1}+o_p(1)
$$
as $d\to \infty$ either when $n$ is fixed or $n\to \infty$.
\end{lemma}
\begin{theorem}
Under (A-i) to (A-iii), it holds that
\begin{equation}
P\Big(\frac{\lambda_1}{\tr(\bSigma)}\in \Big[\frac{(n-1)\tilde{\lambda}_1}{b\tilde{\kappa}+(n-1)\tilde{\lambda}_1},\frac{(n-1)\tilde{\lambda}_1}{a\tilde{\kappa}+(n-1)\tilde{\lambda}_1}\Big]\Big)=1-\alpha+o(1)
\label{2.3}
\end{equation}
as $d\to \infty$ when $n$ is fixed.
\end{theorem}
\begin{remark}
From Theorem 2.1 and Lemma 2.1, under (A-ii) and (A-iii), it holds that $\tr(\bS_D)/\tr(\bSig)=(\tilde{\kappa}+\tilde{\lambda}_1)/\tr(\bSig)=1+o_p(1)$ as $d\to \infty$ and $n\to \infty$.
We have that 
$$
\frac{\tilde{\lambda}_1}{\tr(\bS_D)}=\frac{\lambda_1}{\tr(\bSigma)}\{1+o_p(1)\}.
$$
\end{remark}
\begin{remark}
\label{rem2.3}
The constants $(a,\ b)$ should be chosen for (\ref{2.3}) to have the minimum length. 
If $\lambda_1/\kappa=o(1)$, the length of the confidence interval becomes close to $\{(n-1)\tilde{\lambda}_1/ \tilde{\kappa} \}(1/a-1/b)$ under (A-ii) and (A-iii) when $d\to \infty$ and $n$ is fixed. 
Thus, we recommend to choose constants $(a,\ b)$ such that 
$$
\argmin_{a,b}(1/a-1/b) \quad \mbox{subject to \ $G_{n-1}(b)-G_{n-1}(a)=1-\alpha$},
$$
where $G_{n-1}(\cdot)$ denotes the c.d.f. of $\chi_{n-1}^2$. 
\end{remark}

Let us construct a confidence interval for the contribution ratio of the first principal component.
We used gene expression data by Armstrong et al. (2002) in which the data set consists of $12582\ (=d)$ genes. 
The data set has three leukemia subtypes: 24 samples from acute lymphoblastic leukemia (ALL), 20 samples from mixed-lineage leukemia (MLL), and 28 samples from acute myeloid leukemia (AML). 
We standardized each sample so as to have the unit variance. 
Then, it holds $\tr(\bS)\ (=\tr(\bS_D))=d$, so that $\tilde{\lambda}_{1}+\tilde{\kappa}=d$. 
From Theorem 2.2, we constructed a $95\%$ confidence interval of the first contribution rate for each data set by choosing $(a,\ b)$ as in Remark \ref{rem2.3}. 
The results are summarized in Table 1.
\begin{table}[htb]
{\bf Table 1.} \ The $95\%$ confidence interval (CI) of the first contribution ratio, together with $\tilde{\lambda}_{1}$ and $\tilde{\kappa}$, for Armstrong et al. (2002)'s data sets having $d=12582$.
\begin{center}
\begin{tabular}{c|ccc}
\hline \\[-4mm]
& CI  & $\tilde{\lambda}_{1}$  & $\tilde{\kappa}$ \\
\hline
ALL\ $(n=24)$ & $[0.0557,0.1663]$ &  1256 & 11326\\
MLL\ $(n=20)$ & $[0.1201,0.3458]$ &  2717 & 9865\\
AML\ $(n=28)$ & $[0.0706,0.1884]$ &  1501 & 11081\\
\hline
\end{tabular}
\end{center}
\end{table}
\section{First PC direction and PC score}
In this section, we give asymptotic properties of the first PC direction and PC score in the HDLSS context.
\subsection{Asymptotic properties of the first PC direction}
Let $\hat{\bH}=[\hat{\bh}_1,...,\hat{\bh}_d]$, where $\hat{\bH}$ is a $d\times d$ orthogonal matrix of the sample eigenvectors such that $\hat{\bH}^T\bS\hat{\bH}=\hat{\bLam}$ having $\hat{\bLam}=\mbox{diag}(\hat{\lambda}_1,...,\hat{\lambda}_d)$. 
We assume $\bh_{i}^T\hat{\bh}_{i} \ge 0$ w.p.1 for all $i$ without loss of generality. 
Note that $\hat{\bh}_i$ can be calculated by $\hat{\bh}_i=\{(n-1)\hat{\lambda}_i\}^{-1/2}(\bX-\overline{\bX}) \hat{\bu}_i$.
First, we have the following result. 
\begin{lemma}
Under (A-ii) and (A-iii), it holds that
$$
\hat{\bh}_1^T\bh_1-\Big(1+\frac{\kappa}{\lambda_{1}||\bz_{o1}||^2}\Big)^{-1/2}=o_p(1)
$$
as $d\to \infty$ either when $n$ is fixed or $n\to \infty$.
\end{lemma}

If $\kappa/(n\lambda_1)=o(1)$ as $d\to \infty$ and $n\to \infty$, $\hat{\bh}_1$ is a consistent estimator of $\bh_1$ in the sense that $\hat{\bh}_1^T\bh_1=1+o_p(1)$. 
When $n$ is fixed, $\hat{\bh}_1$ is not a consistent estimator because $\lim_{d\to \infty} \kappa/\lambda_1>0$. 
In order to overcome this inconvenience, we consider applying the NR methodology to the PC direction vector.
Let $\tilde{\bh}_{i}=\{(n-1)\tilde{\lambda}_i\}^{-1/2}(\bX-\overline{\bX})\hat{\bu}_{i}$. 
From Lemma 3.1, we have the following result.
\begin{theorem}
Under (A-ii) and (A-iii), it holds that
$$
\tilde{\bh}_{1}^T\bh_{1}=1+o_p(1)
$$
as $d\to \infty$ either when $n$ is fixed or $n\to \infty$.
\end{theorem}

Note that $||\tilde{\bh}_{1}||^2=\hat{\lambda}_1/\tilde{\lambda}_1\ge 1$ w.p.1.
We emphasize that $\tilde{\bh}_{1}$ is a consistent estimator of $\bh_1$ in the sense of the inner product even when $n$ is fixed though $\tilde{\bh}_{1}$ is not a unit vector. 
We give an application of $\tilde{\bh}_{1}$ in Section 4. 
\subsection{Asymptotic properties of the first PC score}
Let $z_{oij}=z_{ij}-\bar{z}_i$ for all $i,j$. 
First, we have the following result. 
\begin{lemma}
Under (A-ii) and (A-iii), it holds that
$$
\hat{u}_{1j}=z_{o1j}/||\bz_{o1}||+o_p(1)\quad \mbox{for $j=1,...,n$}
$$
as $d\to \infty$ when $n$ is fixed. 
\end{lemma}
\begin{remark}
By using Lemma 3.2 and the test of normality such as Jarque-Bera test, one can check whether (A-i) holds or not.
\end{remark}

By applying the NR methodology to the first PC score, we obtain an estimate by $\tilde{s}_{1j}=\sqrt{(n-1)\tilde{\lambda}_{1}}\hat{u}_{1j},\ j=1,...,n$.
A sample mean squared error of the first PC score is given by MSE$(\tilde{s}_1)=n^{-1}\sum_{j=1}^n(\tilde{s}_{1j}-s_{1j})^2$. 
Then, from Theorem 2.1 and Lemma 3.2, we have the following result.
\begin{theorem}
Under (A-ii) and (A-iii), it holds that
$$
\frac{1}{\sqrt{\lambda_1}}(\tilde{s}_{1j}-s_{1j})=-\bar{z}_{1}+o_p(1) \quad \mbox{for $j=1,...,n$}
$$
as $d\to \infty$ when $n$ is fixed.
Under (A-i) to (A-iii), it holds that
\begin{align*}
\sqrt{\frac{n}{\lambda_1}}(\tilde{s}_{1j}-s_{1j}) \Rightarrow N(0,1) \quad \mbox{for $j=1,...,n$}; \quad \mbox{and} \quad n\frac{\mbox{MSE}(\tilde{s}_1)}{\lambda_1}\Rightarrow \chi_1^2
\end{align*}
as $d\to \infty$ when $n$ is fixed.
\end{theorem}
\begin{remark}
The conventional estimator of the first PC score is given by $\hat{s}_{1j}=\sqrt{(n-1)\hat{\lambda}_{1}}\hat{u}_{1j},\ j=1,...,n$. 
From Theorems 8.1 and 8.2 in Yata and Aoshima (2013), under (A-ii) and (A-iii), it holds that as $d\to \infty$ and $n\to \infty$
\begin{align*}
\frac{\mbox{MSE}(\hat{s}_1)}{\lambda_1}=o_p(1)\ \ \mbox{if $\kappa/(n\lambda_1)=o(1)$},  \quad \mbox{and} \quad 
\frac{\mbox{MSE}(\tilde{s}_1)}{\lambda_1}=o_p(1).
\end{align*}
\end{remark}
\section{Equality tests of two covariance matrices}
In this section, we consider the test of equality of two covariance matrices in the HDLSS context.
Even though there are a variety of tests to deal with covariance matrices when $d\to\infty$ and $n\to \infty$, there seem to be no tests available in the HDLSS context such as $d\to\infty$ while $n$ is fixed. 
Suppose we have two independent $d\times n_{i}$ data matrices, $\bX_{i}=[\bx_{1(i)},...,\bx_{n_i(i)}],\ i=1,2$, where $\bx_{j(i)},\ j= 1,...,n_{i}$, are i.i.d. as a $d$-dimensional distribution, $\pi_i$, having a mean vector $\bmu_i$ and covariance matrix $\bSig_{i}\ (\ge \bO)$.
We assume $n_i \ge 3,\ i=1,2$. 
The eigen-decomposition of $\bSig_i$ is given by $\bSig_i=\bH_i\bLam_i\bH_i^T$, where $\bLam_i=\mbox{diag}(\lambda_{1(i)},...,\lambda_{d(i)})$ having $\lambda_{1(i)}\ge \cdots \ge \lambda_{d(i)}(\ge 0)$ and $\bH_i=[\bh_{1(i)},...,\bh_{d(i)}]$ is an orthogonal matrix of the corresponding eigenvectors. 
\subsection{Equality test using the largest eigenvalues}
We consider the following test for the largest eigenvalues:
\begin{equation}
H_0:\lambda_{1(1)}=\lambda_{1(2)} \quad \mbox{vs.}\quad H_a:\lambda_{1(1)}\neq \lambda_{1(2)}\ \  (\mbox{or }\ H_b:\lambda_{1(1)}< \lambda_{1(2)}). \label{4.1}
\end{equation}
Let $\tilde{\lambda}_{1(i)}$ be the estimate of $\lambda_{1(i)}$ by the NR methodology as in (\ref{2.1}) for $\pi_i$.
Let $\nu_1=n_1-1$ and $\nu_2=n_2-1$. 
From Theorem 2.1, we have the following result. 
\begin{corollary}
Under (A-i) to (A-iii) for each $\pi_i$, it holds that 
$$
\frac{\tilde{\lambda}_{1(1)}/\lambda_{1(1)}}{\tilde{\lambda}_{1(2)}/{\lambda_{1(2)}}} \Rightarrow F_{\nu_1,\nu_2}
$$
as $d\to\infty$ when $n_i$s are fixed, where $F_{\nu_1,\nu_2}$ denotes a random variable distributed as $F$ distribution with degrees of freedom, $\nu_1$ and $\nu_2$.
\end{corollary}
Let $F_1=\tilde{\lambda}_{1(1)}/\tilde{\lambda}_{1(2)}$. 
From Corollary 4.1, we test (\ref{4.1}) for given $\alpha \in(0,1/2)$ by 
\begin{align}
&\mbox{accepting $H_a$}\Longleftrightarrow F_1\notin [ \{F_{\nu_2,\nu_1}(\alpha/2)\}^{-1},F_{\nu_1,\nu_2}(\alpha/2)] \label{4.2}\\
\mbox{or  \quad} &\mbox{accepting $H_b$}\Longleftrightarrow F_1<\{F_{\nu_2,\nu_1}(\alpha)\}^{-1}, \label{4.3}
\end{align}
where $F_{\nu_1,\nu_2}(\alpha)$ denotes the upper $\alpha\%$ point of $F$ distribution with degrees of freedom, $\nu_1$ and $\nu_2$. 
Then, under (A-i) to (A-iii) for each $\pi_i$, it holds that
$$
\mbox{size}=\alpha+o(1)
$$
as $d\to\infty$ when $n_i$s are fixed.

Now, we check the performance of the test by (\ref{4.2}) or (\ref{4.3}).
We also consider a test by the conventional estimator, $\hat{\lambda}_{1(i)}$.
Let $\kappa_{i}=\tr(\bSig_i)-\lambda_{1(i)}=\sum_{s=2}^d\lambda_{s(i)}$ for $i=1,2$. 
From Proposition 2.1, if $\kappa_{i}/\lambda_{1(i)}=o(1)$, $i=1,2$, under (A-i) to (A-iii) for each $\pi_i$ it holds that
$$
\frac{\hat{\lambda}_{1(1)}/\lambda_{1(1)}}{\hat{\lambda}_{1(2)}/{\lambda_{1(2)}}} \Rightarrow F_{\nu_1,\nu_2}
$$
as $d\to\infty$ when $n_i$s are fixed. 
As mentioned in Section 2, the condition `$\kappa_{i}/\lambda_{1(i)}=o(1)$ for $i=1,2$' is quite strict in real high-dimensional data analyses.
Hereafter, we assume $\liminf_{d\to \infty} \kappa_i/\lambda_{1(i)}$
$>0$ for $i=1,2$.  
We analyzed the same gene expression data as in Table 1. 
We set $\alpha=0.05$. 
We considered two cases: (I) $\pi_1:$ ALL ($n_1=24$) and $\pi_2:$ MLL ($n_2=20$), and (II) $\pi_1:$ AML ($n_1=28$) and $\pi_2:$ MLL ($n_2=20$). 
As for $F_1'=\hat{\lambda}_{1(1)}/\hat{\lambda}_{1(2)}$, we considered (\ref{4.2}) and (\ref{4.3}) by replacing $F_1$ with $F_1'$. 
The results are summarized in Table 2.
We observed from Table 2 that only $H_b$ for (I) was accepted by $F_1$, namely, only $F_1$ for (I) rejected $H_0$ vs. $H_b$.
One should note that the condition `$\kappa_{i}/\lambda_{1(i)}=o(1)$ for $i=1,2$' does not hold both for (I) and (II) as observed in Table 1. 
\begin{table}[htb]
{\bf Table 2.} \ Tests of $H_0:\lambda_{1(1)}=\lambda_{1(2)}$ vs. $H_a:\lambda_{1(1)}\neq \lambda_{1(2)}$ or $H_b:\lambda_{1(1)}< \lambda_{1(2)}$ with size $0.05$ for Armstrong et al. (2002)'s data sets having $d=12582$. 
\begin{center}
\begin{tabular}{c|cccc}
\hline
&  $H_a$ by $F_1$ &  $H_a$ by $F_1'$ &  $H_b$ by $F_1$ &  $H_b$ by $F_1'$ \\
\hline
(I) $\pi_1$: ALL, $\pi_2$: MLL  & Reject & Reject & Accept & Reject \\
(II) $\pi_1$: AML, $\pi_2$: MLL & Reject & Reject & Reject & Reject \\
\hline
\end{tabular}
\end{center}
\end{table}
\subsection{Equality test using the largest eigenvalues and their PC directions}
We consider the following test using the largest eigenvalues and their PC directions:
\begin{equation}
H_0:(\lambda_{1(1)},\bh_{1(1)})=(\lambda_{1(2)},\bh_{1(2)}) \quad \mbox{vs.}\quad H_a:(\lambda_{1(1)},\bh_{1(1)})\neq (\lambda_{1(2)},\bh_{1(2)}). 
\label{4.4}
\end{equation}
Let $\tilde{\bh}_{1(i)}$ be the estimator of the first PC direction for $\pi_i$ by the NR methodology given in Section 3.1. 
We assume $\bh_{1(i)}^T\tilde{\bh}_{1(i)} \ge 0$ w.p.1 for $i=1,2$, without loss of generality. 
Here, we have the following result. 
\begin{lemma}
Under (A-ii) and (A-iii) for each $\pi_i$, it holds that 
$$
\tilde{\bh}_{1(1)}^T\tilde{\bh}_{1(2)}=\bh_{1(1)}^T\bh_{1(2)}+o_p(1)
$$
as $d\to\infty$ either when $n_i$ is fixed or $n_i\to \infty$ for $i=1,2$.
\end{lemma}
Let $\tilde{h}=|\tilde{\bh}_{1(1)}^T \tilde{\bh}_{1(2)}|/2+|\tilde{\bh}_{1(1)}^T \tilde{\bh}_{1(2)}|^{-1}/2$. 
Note that $\tilde{h}\ge 1$. 
Then, from Lemma 4.1, we give a test statistic for (\ref{4.4}) as follows:
$$
F_2=\frac{\tilde{\lambda}_{1(1)}}{\tilde{\lambda}_{1(2)}}\tilde{h}_*,
$$
where
$$
\tilde{h}_*=\begin{cases} \tilde{h} & \mbox{if } \tilde{\lambda}_{1(1)} \ge \tilde{\lambda}_{1(2)}, \\
\tilde{h}^{-1} & \mbox{otherwise}.\end{cases}
$$
From Lemma 4.1, we have the following result. 
\begin{theorem}
Under (A-i) to (A-iii) for each $\pi_i$, it holds that
$$
F_2 \Rightarrow F_{\nu_1,\nu_2} \mbox{ under $H_0$}
$$
as $d\to\infty$ when $n_i$s are fixed.
\end{theorem}
From Theorem 4.1, we consider testing (\ref{4.4}) by (\ref{4.2}) with $F_2$ instead of $F_1$.
Then, the size becomes close to $\alpha$ as $d$ increases. 
For the same gene expression data sets as in Section 4.1, we tested (\ref{4.4}) with $\alpha=0.05$ for the cases of (I) and (II).
We observed that only $H_a$ for (II) was accepted by $F_2$, namely, only $F_2$ for (II) rejected $H_0$ vs. $H_a$ in (\ref{4.4}).
\subsection{Equality test of the covariance matrices}
We consider the following test for the covariance matrices:
\begin{equation}
H_0:\bSigma_1=\bSigma_2 \quad \mbox{vs.}\quad H_a:\bSigma_1\neq \bSigma_2. 
\label{4.5}
\end{equation}
When $d\to \infty$ and $n_i$s are fixed, one cannot estimate $\lambda_{j(i)}$s and $\bh_{j(i)}$s for $j=2,...,d$. 
Instead, we consider estimating $\kappa_{i}$s. 
Let $\bS_{D(i)}$ be the dual sample covariance matrix for $\pi_i$.  
We estimate $\kappa_{i}$ by $\tilde{\kappa}_{i}=\tr(\bS_{D(i)})-\tilde{\lambda}_{1(i)}$ for $i=1,2$.
From Lemma 2.1, under (A-ii) and (A-iii) for each $\pi_i$, $\tilde{\kappa}_{i}$s are consistent estimators of $\kappa_{i}$s in the sense that $\tilde{\kappa}_{i}/{\kappa}_{i}=1+o_p(1)$ as $d\to \infty$ when $n_i$s are fixed. 
Let $\tilde{\gamma}=\max\{\tilde{\kappa}_{1}/\tilde{\kappa}_{2}, \tilde{\kappa}_{2}/\tilde{\kappa}_{1}\}$. 
Now, we give a test statistic for (\ref{4.5}) as follows:
$$
F_3=\frac{\tilde{\lambda}_{1(1)}}{\tilde{\lambda}_{1(2)}}\tilde{h}_* \tilde{\gamma }_*,
$$
where
$$
\tilde{\gamma}_*
=\begin{cases} \tilde{\gamma } & \mbox{if } \tilde{\lambda}_{1(1)} \ge \tilde{\lambda}_{1(2)}, \\
\tilde{\gamma }^{-1} & \mbox{otherwise}.\end{cases}
$$
Then, we have the following result. 
\begin{theorem}
Under (A-i) to (A-iii) for each $\pi_i$, it holds that 
$$
F_3 \Rightarrow F_{\nu_1,\nu_2} \mbox{ under $H_0$}
$$
as $d\to\infty$ when $n_i$s are fixed.
\end{theorem}
From Theorem 4.2, we consider testing (\ref{4.5}) by (\ref{4.2}) with $F_3$ instead of $F_1$.
Then, the size becomes close to $\alpha$ as $d$ increases.
For the same gene expression data sets as in Section 4.1, we tested (\ref{4.5}) with $\alpha=0.05$ for the cases of (I) and (II).
We compared the performance of $F_3$ with two other test statistics: $Q_2^2$ and $T_{2}^{2}$ by Srivastava and Yanagihara (2010).
The results are summarized in Table 3.
We observed that $H_a$ was accepted by $F_3$ both for (I) and (II), namely, $F_3$ rejected $H_0$ vs. $H_a$ in (\ref{4.5}) for both the cases.
On the other hand, $Q_2^2$ and $T_2^2$ did not work for these data sets. 
It should be noted that $Q_2^2$ and $T_2^2$ require to meet the conditions that $0<\lim_{d\to\infty} \tr(\bSig^{i})/d <\infty\ (i=1,...,4)$ and $d^{1/2}/n=o(1)$.
As observed in Table 1, the conditions seem not to hold for these data sets with $d=12582$ and $n\le 28$. 
Hence, there is no theoretical guarantee for the results by $Q_2^2$ and $T_2^2$.
\begin{table}[htb]
{\bf Table 3.} \ Tests of $H_0:\bSigma_1=\bSigma_2$ vs. $H_a:\bSigma_1\neq \bSigma_2$ with size $0.05$ for Armstrong et al. (2002)'s data sets having $d=12582$.  
\begin{center}
\begin{tabular}{c|ccc}
\hline \\[-4mm]
& $H_a$ by $F_3$ &  $H_a$ by $Q_2^2$ &  $H_a$ by $T_2^2$  \\
\hline
(I) $\pi_1$: ALL, $\pi_2$: MLL  & Accept & Reject & Reject \\
(II) $\pi_1$: AML, $\pi_2$: MLL & Accept & Reject & Reject \\
\hline
\end{tabular}
\end{center}
\end{table}
\section{Numerical results and discussions}
\subsection{Comparisons of the estimates on the first PC}
In this section, we compared the performance of $\tilde{\lambda}_1$, $\tilde{\bh}_1$ and $\tilde{s}_{1j}$ with their conventional counterparts by Monte Carlo simulations. 
We set $d=2^{k},\ k=3,...,11$ and $n=10$. 
We considered two cases for $\lambda_i$s: (a) $\lambda_i=d^{1/i}$, $i=1,...,d$ and (b) $\lambda_i=d^{3/(2+2i)}$, $i=1,...,d$. 
Note that $\lambda_1=d$ for (a) and $\lambda_1=d^{3/4}$ for (b). 
Also, note that (A-ii) holds both for (a) and (b). 
Let $d_*= \lceil d^{1/2} \rceil $, where $\lceil x \rceil$ denotes the smallest integer $\ge x$. 
We considered a non-Gaussian distribution as follows: $(z_{1j},...,z_{d-d_* j})^T,\ j=1,...,n,$ are i.i.d. as $N_{d-d_*}(\bze, \bI_{d-d_*})$ and $(z_{d-d_*+1 j},...,z_{d j})^T,\ j=1,...,n,$ are i.i.d. as the $d_*$-variate $t$-distribution, $t_{d_*}(\bze,\bI_{d_*},10)$ with mean zero, covariance matrix $\bI_{d_*}$ and degrees of freedom $10$, where $(z_{1j},...,z_{d-d_* j})^T$ and $(z_{d-d_*+1j},...,z_{dj})^T$ are independent for each $j$. 
Note that (A-i) and (A-iii) hold both for (a) and (b) from the fact that $\sum_{r,s\ge 2}^d \lambda_r\lambda_s E\{(z_{rk}^2-1)(z_{sk}^2-1)\}=2\sum_{s= 2}^{d-d_*}\lambda_s^2+O(\sum_{r,s\ge d-d_*+1}^d \lambda_r\lambda_{s})=o(\lambda_1^2)$. 

The findings were obtained by averaging the outcomes from $2000\ (=R$, say) replications. 
Under a fixed scenario, suppose that the $r$-th replication ends with estimates, ($\hat{\lambda}_{1r}$, $\hat{\bh}_{1r}$, MSE$(\hat{s}_1)_r$) and ($\tilde{\lambda}_{1r}$, $\tilde{\bh}_{1r}$, MSE$(\tilde{s}_1)_r$) $(r=1,...,R)$.
Let us simply write $\hat{\lambda}_1=R^{-1}\sum_{r=1}^R \hat{\lambda}_{1r}$ and $\tilde{\lambda}_1=R^{-1}\sum_{r=1}^R \tilde{\lambda}_{1r}$. 
We also considered the Monte Carlo variability by $\mbox{var}(\hat{\lambda}_1/\lambda_1)=(R-1)^{-1}\sum_{r=1}^R(\hat{\lambda}_{1r}-\hat{\lambda}_1)^2/\lambda_1^{2}$ and 
$\mbox{var}(\tilde{\lambda}_1/\lambda_1)=(R-1)^{-1}\sum_{r=1}^R(\tilde{\lambda}_{1r}-\tilde{\lambda}_1)^2/\lambda_1^{2}$. 
Figure 1 shows the behaviors of ($\hat{\lambda}_1/\lambda_1$, $\tilde{\lambda}_1/\lambda_1$) in the left panel and (var$(\hat{\lambda}_1/\lambda_1)$, var$(\tilde{\lambda}_1/\lambda_1)$) in the right panel for (a) and (b). 
We gave the asymptotic variance of $\tilde{\lambda}_1/\lambda_1$ by Var$\{\chi_{n-1}^2/(n-1)\}=0.222$ from Theorem 2.1 and showed it by the solid line in the right panel. 
We observed that the sample mean and variance of $\tilde{\lambda}_1/\lambda_1$ become close to those asymptotic values as $d$ increases. 
\begin{figure}
\includegraphics[scale=0.52]{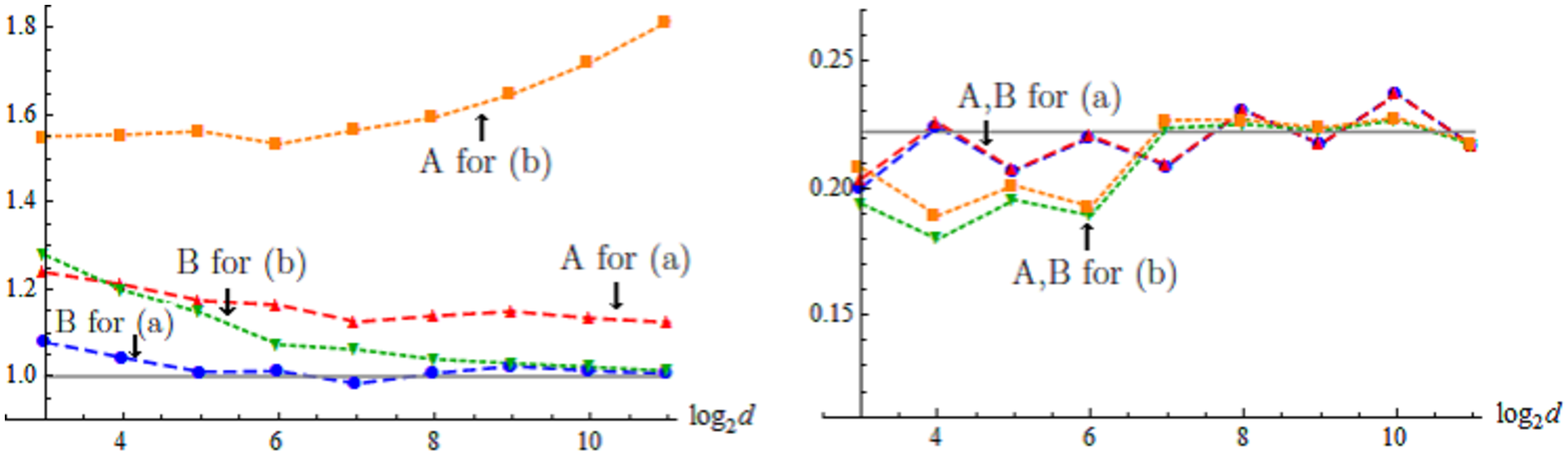} \\
{\small \hspace{1.5cm} A: $\hat{\lambda}_1/\lambda_1$ and B: $\tilde{\lambda}_1/\lambda_1$ \hspace{2.2cm}  A: var$(\hat{\lambda}_1/\lambda_1)$ and B: var$(\tilde{\lambda}_1/\lambda_1)$}
\caption{The values of A: $\hat{\lambda}_1/\lambda_1$ and B: $\tilde{\lambda}_1/\lambda_1$ are denoted by the dashed lines for (a) and by the dotted lines for (b) in the left panel.
The values of A: var$(\hat{\lambda}_1/\lambda_1)$ and B: var$(\tilde{\lambda}_1/\lambda_1)$ are denoted by the dashed lines for (a) and by the dotted lines for (b) in the left panel. 
The asymptotic variance of $\tilde{\lambda}_1/\lambda_1$ was given by Var$\{\chi_{n-1}^2/(n-1)\}=0.222$ and denoted by the solid line in the left panel. 
}
\end{figure}

Similarly, we plotted ($\hat{\bh}_1^T\bh_1$, $\tilde{\bh}_1^T\bh_1$) and (var$(\hat{\bh}_1^T\bh_1)$, var$(\tilde{\bh}_1^T\bh_1)$) in Figure 2 and (MSE$(\hat{s}_1)/\lambda_1$, MSE$(\tilde{s}_1)/\lambda_1$) and (var(MSE$(\hat{s}_1)/\lambda_1$), var(MSE$(\tilde{s}_1)/\lambda_1$)) in Figure 3.
From Theorem 3.2, we gave the asymptotic mean of MSE$(\tilde{s}_1)/\lambda_1$ by $E(\chi_{1}^2/n)=0.1$ and showed it by the solid line in the left panel of Figure 3.
We also gave the asymptotic variance of MSE$(\tilde{s}_1)/\lambda_1$ by Var$(\chi_{1}^2/n)=0.02$ in the right panel of Figure 3. 
Throughout, the estimators by the NR method gave good performances both for (a) and (b) when $d$ is large. 
However, the conventional estimators gave poor performances especially for (b). 
This is probably because the bias of the conventional estimators, $\kappa/(n\lambda_1)$, is large for (b) compared to (a). 
See Proposition 2.1 for the details. 
\begin{figure}
\includegraphics[scale=0.52]{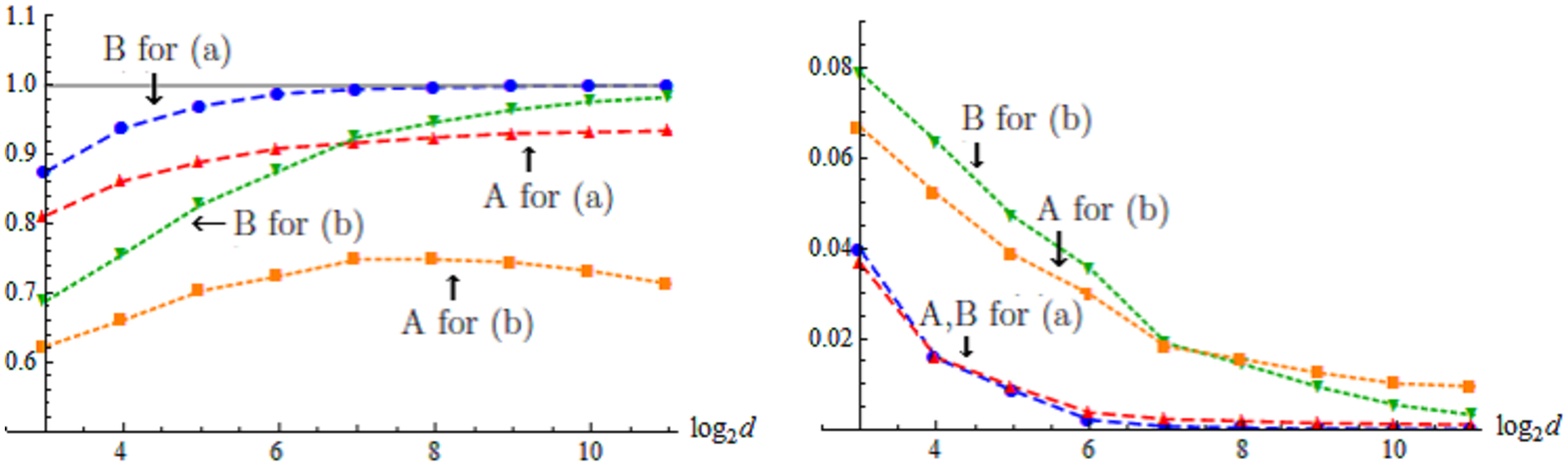} \\
{\small
\hspace{1.5cm} A: $\hat{\bh}_1^T\bh_1$ and B: $\tilde{\bh}_1^T\bh_1$ \hspace{2.4cm} A: var$(\hat{\bh}_1^T\bh_1)$ and B: var$(\tilde{\bh}_1^T\bh_1)$}
\caption{The values of A: $\hat{\bh}_1^T\bh_1$ and B: $\tilde{\bh}_1^T\bh_1$ are denoted by the dashed lines for (a) and by the dotted lines for (b) in the left panel.
The values of A: var$(\hat{\bh}_1^T\bh_1)$ and B: var$(\tilde{\bh}_1^T\bh_1)$ are denoted by the dashed lines for (a) and by the dotted lines for (b) in the right panel. 
}
\vspace{3mm}
\includegraphics[scale=0.52]{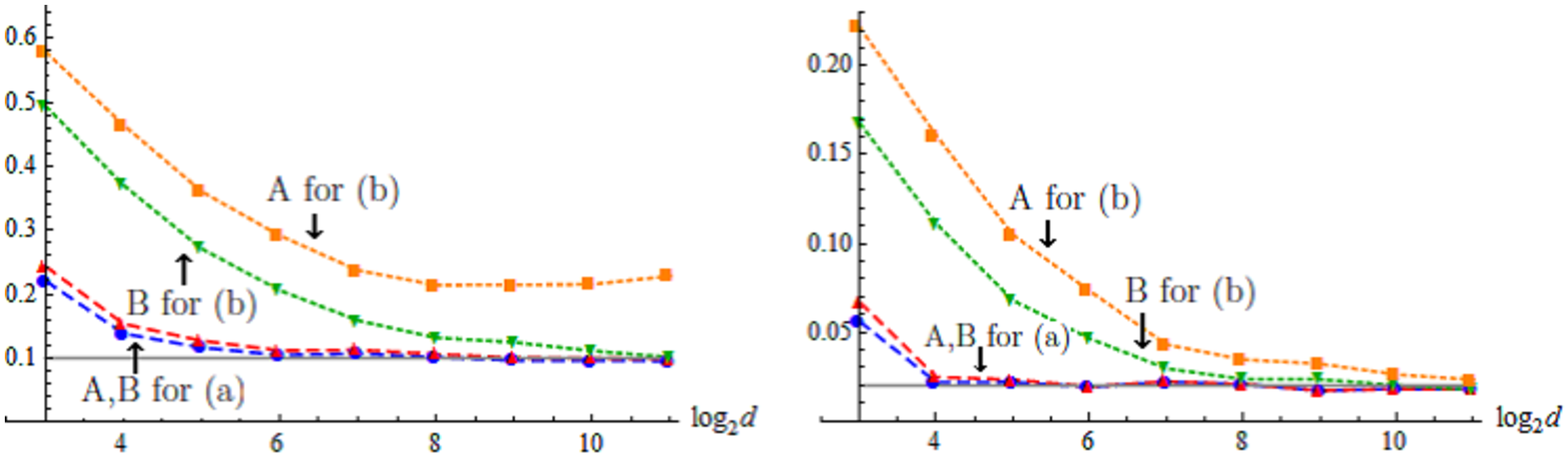} \\
{\small
A: MSE$(\hat{s}_1)/\lambda_1$ and B: MSE$(\tilde{s}_1)/\lambda_1$ \quad \ \ \ A: var(MSE$(\hat{s}_1)/\lambda_1$) and B: var(MSE$(\tilde{s}_1)/\lambda_1$)}
\caption{The values of A: MSE$(\hat{s}_1)/\lambda_1$ and B: MSE$(\tilde{s}_1)/\lambda_1$ are denoted by the dashed lines for (a) and by the dotted lines for (b) in the left panel.
The values of A: var(MSE$(\hat{s}_1)/\lambda_1$) and B: var(MSE$(\tilde{s}_1)/\lambda_1$) are denoted by the dashed lines for (a) and by the dotted lines for (b) in the right panel.
The asymptotic mean and variance of MSE$(\tilde{s}_1)/\lambda_1$ were given by $E(\chi_{1}^2/n)=0.1$ and Var$(\chi_{1}^2/n)=0.02$ and denoted by the solid lines in both the panels.
}
\end{figure}
\subsection{Equality tests of two covariance matrices}
We used computer simulations to study the performance of the test procedures by $F_1$ for (\ref{4.1}), $F_2$ for (\ref{4.4}) and $F_3$ for (\ref{4.5}). 
We set $\alpha=0.05$. 
Independent pseudo-random normal observations were generated from $\pi_i:N_d(\bze,\bSig_i)$,\ $i=1,2$. 
We set $(n_1,n_2)=(10,20)$.
We considered the cases: $d=2^k,\ k=3,...,11$, and 
\begin{equation}
\bSig_{i}=\left( \begin{array}{cc}
\bSig_{i(1)} & \bO_{2, d-2} \\
\bO_{d-2, 2} & \bSig_{i(2)}
\end{array} \right),\ i=1,2,\label{5.1}
\end{equation}
where $\bO_{k,l}$ is the $k\times l$ zero matrix, $\bSig_{1(1)}=\mbox{diag}(d^{3/4},d^{1/2})$ and $\bSig_{1(2)}=(0.3^{|s-t|})$.
When considered the alternative hypotheses, we set 
\begin{equation}
\bSig_{2(1)}=\left( \begin{array}{cc}
1/3 & \sqrt{8}/3  \\
\sqrt{8}/3  & -1/3
\end{array} \right)\mbox{diag}(3d^{3/4},1.5d^{1/2})\left( \begin{array}{cc}
1/3 & \sqrt{8}/3  \\
\sqrt{8}/3  & -1/3
\end{array} \right)
\label{5.2}
\end{equation}
and $\bSig_{2(2)}=1.5 (0.3^{|s-t|})$.
Note that $\lambda_{1(2)}/\lambda_{1(1)}=3$, $\kappa_2/\kappa_1=1.5$, $\bh_{1(1)}=(1,0,....,0)^T$ and $\bh_{1(2)}=(1/3,\sqrt{8}/3,0....,0)^T$, so that $\bh_{1(1)}^T\bh_{1(2)}=1/3$. 
Also, note that (A-i) to (A-iii) hold for each $\pi_i$. 
Let $h=(|\bh_{1(1)}^{T}\bh_{1(2)}|+1/|\bh_{1(1)}^{T}\bh_{1(2)}|)/2$ and $\gamma=\max\{\kappa_{1}/\kappa_{2},\kappa_{2}/\kappa_{1}\}$. 
From Lemmas 2.1 and 4.1, it holds that $\tilde{h}=h+o_p(1)$ and $\tilde{\gamma}=\gamma+o_p(1)$. 
Thus, from Corollary 4.1, Theorems 4.1 and 4.2, we obtained the asymptotic powers of $F_1$, $F_2$ and $F_3$ with $(\tilde{h}_{*}, \tilde{\gamma}_{*})=(h^{-1}, \gamma^{-1})$ as follows:
\begin{align*}
&\mbox{Power}(F_1)=P\big\{(\lambda_{1(1)}/\lambda_{1(2)} )f \notin [\{F_{\nu_2,\nu_1}(\alpha/2)\}^{-1},F_{\nu_1,\nu_2}(\alpha/2)] \big\}=0.39, \\
&\mbox{Power}(F_2)=P\big\{h^{-1}(\lambda_{1(1)}/\lambda_{1(2)} )f \notin [\{F_{\nu_2,\nu_1}(\alpha/2)\}^{-1},F_{\nu_1,\nu_2}(\alpha/2)] \big\}=0.726 \\
&\mbox{and \ Power}(F_3)=P\big\{\gamma^{-1}h^{-1}(\lambda_{1(1)}/\lambda_{1(2)} )f \notin [\{F_{\nu_2,\nu_1}(\alpha/2)\}^{-1},F_{\nu_1,\nu_2}(\alpha/2)] \big\}\\
&\qquad \qquad \qquad \ =0.908,
\end{align*}
where $f$ denotes a random variable distributed as $F$ distribution with degrees of freedom, $\nu_1$ and $\nu_2$.
Note that $\mbox{Power}(F_2)$ and $\mbox{Power}(F_3)$ give lower bounds of the asymptotic powers when $\tilde{h}_{*}=h^{-1}$ and $\tilde{\gamma}_{*}=\gamma^{-1}$. 

In Figure 4, we summarized the findings obtained by averaging the outcomes from 4000 $(=R,$ say) replications.
Here, the first $2000$ replications were generated by setting $\bSig_2=\bSig_1$ as in (\ref{5.1}) and the last $2000$ replications were generated by setting $\bSig_2$ as in (\ref{5.2}). 
Let $F_{ir}\ (i=1,2,3)$ be the $r$th observation of $F_{i}$ for $r=1,...,4000$. 
We defined $P_{r}=1\ (\mbox{or}\ 0)$ when $H_0$ was falsely rejected (or not) for $r=1,...,2000$, and $H_a$ was falsely rejected (or not) for $r=2001,...,4000$.  
We defined $\overline{\alpha}=(R/2)^{-1} \sum_{r=1}^{R/2}P_{r}$ to estimate the size and $1-\overline{\beta}=1-(R/2)^{-1} \sum_{r=R/2+1}^{R}P_{r}$ to estimate the power. 
Their standard deviations are less than $0.011$. 
Throughout, the tests gave adequate performances for the high-dimensional cases. 
\begin{figure}
\includegraphics[scale=0.52]{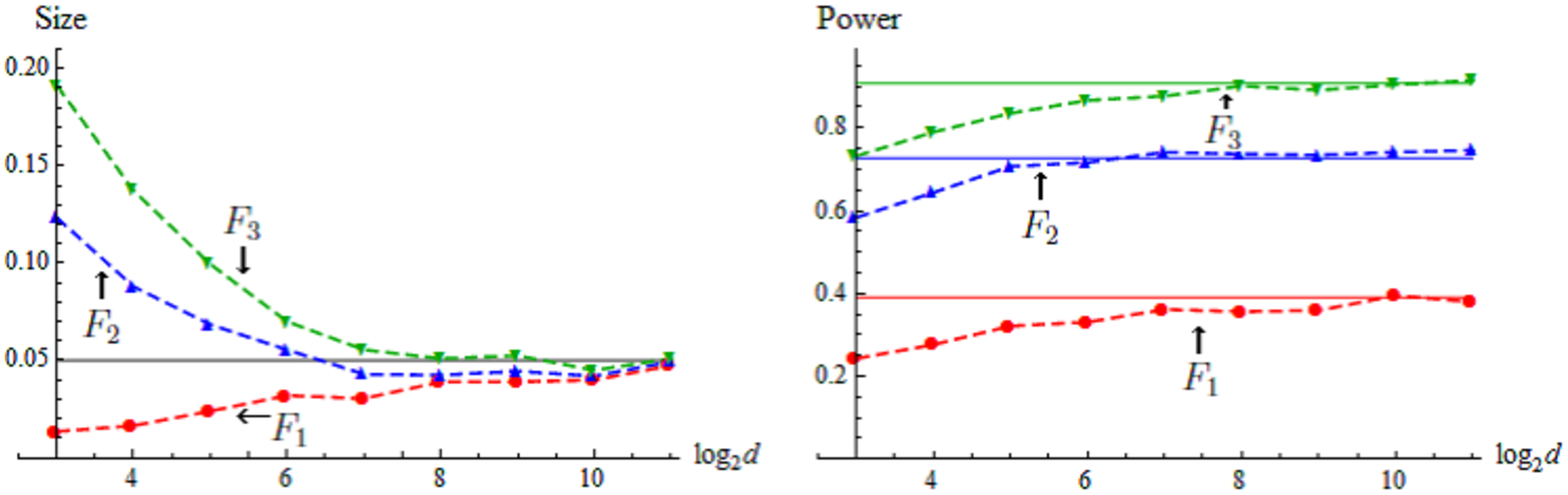} \\
{\small \hspace{1.5cm} Sizes of $F_1$, $F_2$ and $F_3$ \hspace{2.6cm} 
Powers of $F_1$, $F_2$ and $F_3$}
\caption{The values of $\overline{\alpha}$ are denoted by the dashed lines in the left panel and the values of $1-\overline{\beta}$ are denoted by the dashed lines in the right panel for $F_1$, $F_2$ and $F_3$. 
The asymptotic powers were given by $\mbox{Power}(F_1)=0.39$, $\mbox{Power}(F_1)=0.726$ and $\mbox{Power}(F_3)=0.908$ which were denoted by the solid lines in the right panel. 
}
\end{figure}
\appendix
\section{}
Throughout, let $\bP_n=\bI_n-\bone_n\bone_n^T/n$, where $\bone_n=(1,...,1)^T$.
Let $\be_{n}=(e_{1},...,e_{n})^T$ be an arbitrary (random) $n$-vector such that $||\be_{n} ||=1$ and $\be_{n}^T \bone_n=0$.
\\[5mm]
{\it Proof of Proposition 2.1.} \ 
We assume $\bmu=\bze$ without loss of generality. 
We write that $\bX^T\bX=\sum_{s=1}^{i_*}\lambda_s \bz_{s}\bz_{s}^T+\sum_{s=i_*+1}^d\lambda_s \bz_{s}\bz_{s}^T$ for $i_*=1$ when $n$ is fixed, and for some fixed $i_*(\ge1)$ when $n\to\infty$.  
Here, by using Markov's inequality, for any $\tau>0$, under (A-ii) and (A-iii), we have that 
\begin{align*}
&P\Big\{ \sum_{j=1}^{n} \Big(  \sum_{s=i_*+1}^d \frac{\lambda_s (z_{sj}^2-1)}{n \lambda_1 } \Big)^2>\tau \Big\} \le  \frac{\sum_{r,s\ge 2}^d \lambda_r\lambda_{s}E\{(z_{rk}^2-1)(z_{sk}^2-1)\}} { \tau n\lambda_1^2} \to 0\\
&\mbox{and }\ P\Big\{ \sum_{j\neq j'}^{n} \Big(  \sum_{s=i_*+1}^d \frac{\lambda_s z_{sj}z_{sj'}}{n \lambda_1 } \Big)^2>\tau \Big\} 
\le  \frac{\delta_{i_*}}{\tau \lambda_1^2} \to 0
\end{align*}
as $d\to \infty$ either when $n$ is fixed or $n\to \infty$. 
Note that $\sum_{j=1}^ne_{j}^4\le 1$ and $\sum_{j\neq j'}^ne_{j}^2e_{j'}^2\le 1$.
Then, under (A-ii) and (A-iii), we have that 
\begin{align*}
\Big| \sum_{j=1}^ne_{j}^2\sum_{s=i_*+1}^d \frac{ \lambda_s (z_{sj}^2-1) }{n \lambda_1 }
\Big|&\le \Big\{\sum_{j=1}^ne_{j}^4 \Big\}^{1/2}   \Big\{ \sum_{j=1}^n \Big(  \sum_{s=i_*+1}^d \frac{ \lambda_s (z_{sj}^2-1) }{n \lambda_1 } \Big)^2\Big\}^{1/2} \\
&=o_p(1) \quad \mbox{and} \\
\Big|\sum_{j\neq j'}^{n}e_je_{j'}  \sum_{s=i_*+1}^d \frac{ \lambda_s z_{sj}z_{sj'} }{n \lambda_1 } \Big|&\le 
\Big\{\sum_{j\neq j'}^ne_{j}^2e_{j'}^2 \Big\}^{1/2} \Big\{\sum_{j\neq j'}^{n} \Big(  \sum_{s=i_*+1}^d \frac{ \lambda_s z_{sj}z_{sj'} }{n \lambda_1 } \Big)^2
 \Big\}^{1/2} \\
&=o_p(1)
\end{align*}
as $d\to \infty$ either when $n$ is fixed or $n\to \infty$.
Thus, we claim that 
\begin{equation}
\be_{n}^T\frac{\bX^T\bX}{(n-1)\lambda_1}\be_{n}=\be_{n}^T\frac{\sum_{s=1}^{i_*}\lambda_s \bz_{s}\bz_{s}^T}{(n-1)\lambda_1}\be_{n}
+\frac{\kappa}{(n-1)\lambda_1}+o_p(1) \label{A.1}
\end{equation}
from the fact that $\sum_{s=i_*+1}^{d}\lambda_s/(n\lambda_1)=\kappa/(n\lambda_1)+o(1)$ when $n\to \infty$. 
Note that $\be_{n}^T\bP_n=\be_{n}^T$ and $ \bP_n \bz_{s}=\bz_{os}$ for all $s$. 
Also, note that $(\bz_{os}/n^{1/2})^T(\bz_{os'}/n^{1/2})=o_p(1)$ for $s\neq s'$ as $n\to \infty$ from the fact that $E\{(\bz_{os}^T\bz_{os'}/n)^2\}=o(1)$ as $n\to \infty$.
Then, by noting that $P(\lim_{d\to \infty} ||\bz_{o1}|| \neq 0 )=1$, $\liminf_{d\to \infty}\lambda_{1}/\lambda_{2}>1$ and $\bz_{o1}^T\bone_n=0$, it holds that 
\begin{align}
\max_{\be_{n}}\Big\{
\be_{n}^T\frac{\sum_{s=1}^{i_*}\lambda_s \bz_{s}\bz_{s}^T}{(n-1)\lambda_1}\be_{n}\Big\}&=
\max_{\be_{n}}\Big\{ \be_{n}^T\frac{\sum_{s=1}^{i_*}\lambda_s \bz_{os}\bz_{os}^T}{(n-1)\lambda_1}\be_{n}\Big\} \notag \\
&=||\bz_{o1}/\sqrt{n-1}||^2+o_p(1) 
\label{A.2}
\end{align}
as $d\to \infty$ either when $n$ is fixed or $n\to \infty$. 
Note that $\hat{\bu}_1^T \bone_n=0$ and $\hat{\bu}_1^T\bP_n=\hat{\bu}_1^T$ when $\bS_{D} \neq \bO$. 
Then, from (\ref{A.1}), (\ref{A.2}) and $\bP_n\bX^T\bX \bP_n/(n-1)=\bS_{D}$, under (A-ii) and (A-iii), we have that 
\begin{equation}
\hat{\bu}_1^T\frac{\bS_{D}}{\lambda_1}\hat{\bu}_1=\hat{\bu}_1^T\frac{\bX^T\bX}{(n-1)\lambda_1}\hat{\bu}_1
=||\bz_{o1}/\sqrt{n-1}||^2
+\frac{\kappa}{(n-1) \lambda_1}+o_p(1)\label{A.3}
\end{equation}
as $d\to \infty$ either when $n$ is fixed or $n\to \infty$.
It concludes the result. 
\hfill $\Box$
\\[5mm]
{\it Proof of Lemma 2.1.} \ 
By using Markov's inequality, for any $\tau>0$, under (A-ii) and (A-iii), we have that 
\begin{align*}
&P\Big\{  \Big(  \sum_{s=2}^d \frac{\lambda_s \{||\bz_{os}||^2-(n-1)\} }{ (n-1) \lambda_1 } \Big)^2>\tau \Big\} \\
&=P\Big\{  \Big(  \sum_{s=2}^d \frac{\lambda_s \{(n-1)\sum_{k=1}^n(z_{sk}^2-1)/n-\sum_{k\neq k'}^nz_{sk}z_{sk'}/n \} }{ (n-1) \lambda_1 } \Big)^2>\tau \Big\} 
\\
&=O\Big\{ \frac{\sum_{r,s\ge 2}^d \lambda_r\lambda_{s}E\{(z_{rk}^2-1)(z_{sk}^2-1)\}} {  n\lambda_1^2}  \Big\}+O\{\delta_1/(n\lambda_1)^2\}\to 0 
\end{align*}
as $d\to \infty$ either when $n$ is fixed or $n\to \infty$. 
Thus it holds that $\tr(\bS_{D})/\lambda_1=\kappa/\lambda_1+  ||\bz_{o1}/\sqrt{n-1}||^2+o_p(1)$ from the fact that $\tr(\bS_{D})=\lambda_1 ||\bz_{o1}||^2/(n-1)+ \sum_{s=2}^d \lambda_s ||\bz_{os}||^2/(n-1)$. 
Then, from Proposition 2.1 and $\liminf_{d\to \infty}\kappa/\lambda_1>0$, we can claim the results. 
\hfill $\Box$ 
\\[5mm]
{\it Proof of Theorem 2.1.} \ 
When $n\to \infty$, we can claim the results from Theorems 4.1, 4.2 and Corollary 4.1 in Yata and Aoshima (2013).
When $n$ is fixed, by combining Proposition 2.1 with Lemma 2.1, we can claim the results because $||\bz_{o1}||^2=\sum_{k=1}^nz_{1k}^2-n\bar{z}_1^2$ is distributed as $\chi_{n-1}^2$ if $z_{1j},\ j=1,...,k,$ are i.i.d. as $N(0,1)$.
\hfill $\Box$ 
\\[5mm]
{\it Proof of Theorem 2.2.} \
From Theorem 2.1 and Lemma 2.1, under (A-i) to (A-iii), it holds that 
\begin{align*}
&P\Big(\frac{\lambda_1}{\tr(\bSigma)}\in \Big[\frac{(n-1)\tilde{\lambda}_1}{b\tilde{\kappa}+(n-1)\tilde{\lambda}_1},\frac{(n-1)\tilde{\lambda}_1}{a\tilde{\kappa}+(n-1)\tilde{\lambda}_1}\Big]\Big)\\
&=
P\Big(\frac{(n-1)\tilde{\lambda}_1}{b \tilde{\kappa}+(n-1)\tilde{\lambda}_1}\le \frac{\lambda_1}{\tr(\bSigma)} \le \frac{(n-1)\tilde{\lambda}_1}{a \tilde{\kappa}+(n-1)\tilde{\lambda}_1} \Big)  \\
&=P\Big(\frac{a\tilde{\kappa}}{ (n-1) \tilde{\lambda}_1}  \le \frac{\kappa}{\lambda_1}\le 
\frac{b\tilde{\kappa}}{(n-1)\tilde{\lambda}_1} \Big)=P\Big(a\le (n-1) \frac{\tilde{\lambda}_1 \kappa}{\lambda_1 \tilde{\kappa}} \le b \Big)\\
&=1-\alpha+o(1)
\end{align*}
as $d\to \infty$ when $n$ is fixed. 
It concludes the result.
\hfill $\Box$ 
\\[5mm]
{\it Proofs of Lemmas 3.1 and 3.2.} \ 
We note that $||\bz_{o1}||^2/n=1+o_p(1)$ as $n\to \infty$. 
From (\ref{A.3}), under (A-ii) and (A-iii), we have that 
\begin{equation}
\hat{\bu}_{1}^T\bz_{o1}/||\bz_{o1}||=1+o_p(1) \label{A.4}
\end{equation}
as $d\to \infty$ either when $n$ is fixed or $n\to \infty$, so that $\hat{\bu}_{1}^T\bz_{o1}=||\bz_{o1}||+o_p(n^{1/2})$.
Thus, we can claim the result of Lemma 3.2. 
On the other hand, with the help of Proposition 2.1, under (A-ii) and (A-iii), it holds that from (\ref{A.4})
\begin{align*}
\bh_1^T \hat{\bh}_1=
\frac{\bh_1^T (\bX-\overline{\bX}) \hat{\bu}_1}{\{(n-1)\hat{\lambda}_1\}^{1/2}}=
\frac{\lambda_1^{1/2}\bz_{o1}^T \hat{\bu}_1}{\{(n-1)\hat{\lambda}_1\}^{1/2}}
&=\frac{||\bz_{o1}||+o_p(n^{1/2})}{\{||{\bz}_{o1}||^2+\kappa/\lambda_1+o_p(n)\}^{1/2}}\\
&=\frac{1}{\{1+\kappa/(\lambda_1||{\bz}_{o1}||^2)\}^{1/2}}+o_p(1)
\end{align*}
as $d\to \infty$ either when $n$ is fixed or $n\to \infty$.
It concludes the result of Lemma 3.1.
\hfill $\Box$ 
\\[5mm]
{\it Proof of Theorem 3.1.} \ 
With the help of Theorem 2.1, under (A-ii) and (A-iii), we have that from (\ref{A.4})
$$
\bh_1^T \tilde{\bh}_1=
\frac{\bh_1^T (\bX-\overline{\bX}) \hat{\bu}_1}{\{(n-1)\tilde{\lambda}_1\}^{1/2}}
=\frac{||\bz_{o1}||+o_p(n^{1/2})}{\{||{\bz}_{o1}||^2+o_p(n)\}^{1/2}}
=1+o_p(1)
$$
as $d\to \infty$ either when $n$ is fixed or $n\to \infty$.
It concludes the result. 
\hfill $\Box$
\\[5mm]
{\it Proof of Theorem 3.2.} \ 
By combing Theorem 2.1 with Lemma 3.2, under (A-ii) and (A-iii), we have that 
$$
\tilde{s}_{1j}/\sqrt{\lambda_1}=\hat{u}_{1j}\sqrt{(n-1)\tilde{\lambda}_1/\lambda_1}=\hat{u}_{1j}||\bz_{o1}||+o_p(1)=z_{o1j}+o_p(1)
$$
as $d\to \infty$ when $n$ is fixed. 
By noting that $z_{o1j}=z_{1j}-\bar{z}_1$ and $\bar{z}_1$ is distributed as $N(0,1/n)$ under (A-i), we have the results. 
\hfill $\Box$
\\[5mm]
{\it Proof of Corollary 4.1}. \ 
From Theorem 2.1, the result is obtained straightforwardly. 
\hfill $\Box$
\\[5mm]
{\it Proof of Lemma 4.1.} \ 
Let $\bZ_{i}=[\bz_{1(i)},...,\bz_{d(i)}]^T$ be a sphered data matrix of $\pi_i$ for $i=1,2$, where $\bz_{j(i)}=(z_{j1(i)},...,z_{jn_i(i)})^T$. 
We assume $\bmu_1=\bmu_2=\bze$ without loss of generality. 
Let $\beta_{st}=( \lambda_{s(1)}\lambda_{t(2)})^{1/2} \bh_{s(1)}^T\bh_{t(2)}$ for all $s,t$. 
Let $i_{\star}$ be a fixed constant such that $\sum_{s=i_{\star}+1}^d\lambda_{s(j)}^2/\lambda_{1(j)}^2=o(1)$ as $d\to \infty$ for $j=1,2$.
Note that $i_{\star}$ exists under (A-ii) for each $\pi_i$.
We write that 
\begin{align*}
\bX_1^T\bX_2&=\sum_{s,t\le i_{\star}}\beta_{st}\bz_{s(1)}\bz_{t(2)}^T +
\sum_{s,t\ge i_{\star}+1}^d\beta_{st} \bz_{s(1)}\bz_{t(2)}^T\\
&\quad + \sum_{s= i_{\star}+1}^d\sum_{t=1}^{i_{\star}} \beta_{st} \bz_{s(1)}\bz_{t(2)}^T
+\sum_{s=1}^{i_{\star}} \sum_{t=i_{\star}+1}^d \beta_{st} \bz_{s(1)}\bz_{t(2)}^T.
\end{align*}
Note that 
\begin{align*}
&E\Big\{ \Big( \sum_{s=i_{\star}+1}^d\sum_{t=1}^{i_{\star}} \beta_{st} z_{sj(1)}z_{tj'(2)} \Big)^2\Big\} \\
&=\tr\Big(\sum_{s=i_{\star}+1}^d \lambda_{s(1)}\bh_{s(1)}\bh_{s(1)}^T\sum_{t=1}^{i_{\star}} \lambda_{t(2)}\bh_{t(2)}\bh_{t(2)}^T \Big)
\le i_{\star} \lambda_{i_{\star}+1(1)}\lambda_{1(2)}
\end{align*}
for all $j,j'$. 
Also, note that 
\begin{align*}
E\Big\{ \Big( \sum_{s,t\ge i_{\star}+1}^d  \beta_{st}  z_{sj(1)}z_{tj'(2)} \Big)^2\Big\} 
&=\tr\Big(\sum_{s=i_{\star}+1}^d \lambda_{s(1)}\bh_{s(1)}\bh_{s(1)}^T\sum_{t=i_{\star}+1}^d \lambda_{t(2)}\bh_{t(2)}\bh_{t(2)}^T \Big) \\
&\le 
\Big( \sum_{s=i_{\star}+1}^d\lambda_{s(1)}^2\sum_{t=i_{\star}+1}^d\lambda_{t(2)}^2 \Big)^{1/2}
\end{align*}
for all $j,j'$. 
Then, by using Markov's inequality, for any $\tau>0$, under (A-ii) for each $\pi_i$, we have that 
\begin{align*}
&P\Big\{ \sum_{j=1}^{n_1}\sum_{j'=1}^{n_2} 
\Big( \sum_{s=i_{\star}+1}^d\sum_{t=1}^{i_{\star}}  \frac{ \beta_{st} z_{sj(1)}z_{tj'(2)}}
{ ( n_1 n_2 \lambda_{1(1)}\lambda_{1(2)})^{1/2} } \Big)^2>\tau \Big\} \to 0, \\
&P\Big\{ \sum_{j=1}^{n_1}\sum_{j'=1}^{n_2} 
\Big( \sum_{s=1}^{i_{\star}}\sum_{t=i_{\star}+1}^{d} \frac{ \beta_{st}  z_{sj(1)}z_{tj'(2)}}
{ ( n_1 n_2 \lambda_{1(1)}\lambda_{1(2)})^{1/2} } \Big)^2>\tau \Big\} \to 0 \\
&\mbox{and } P\Big\{ \sum_{j=1}^{n_1}\sum_{j'=1}^{n_2} 
\Big(  \sum_{s,t\ge i_{\star}+1}^d \frac{ \beta_{st} z_{sj(1)}z_{tj'(2)}}
{ ( n_1 n_2 \lambda_{1(1)}\lambda_{1(2)})^{1/2} } \Big)^2>\tau \Big\} \to 0
\end{align*}
as $d\to \infty$ either when $n_i$ is fixed or $n_i\to \infty$ for $i=1,2$. 
Hence, similar to (\ref{A.1}), it holds that 
$$
\frac{\be_{n_1}^T \bX_1^T\bX_2 \be_{n_2} }{( \nu_1 \nu_2 \lambda_{1(1)}\lambda_{1(2)})^{1/2}}
=
\frac{ \be_{n_1}^T \sum_{s,t \le i_{\star}} \beta_{st} \bz_{s(1)}\bz_{t(2)}^T\be_{n_2}}{( \nu_1 \nu_2 \lambda_{1(1)}\lambda_{1(2)})^{1/2}}+o_p(1).
$$
Note that $\be_{n_i}^T\bP_{n_i}=\be_{n_i}^T$ and $ \bP_{n_i} \bz_{1(i)}=\bz_{o1(i)}$ for $i=1,2$, 
where $\bz_{o1(i)}=\bz_{1(i)}-(\bar{z}_{1(i)},...,\bar{z}_{1(i)})^T$ and $\bar{z}_{1(i)}=n_{i}^{-1}\sum_{k=1}^{n_i}z_{1k(i)}$. 
Also, note that $\bX_i\bP_{n_i}=(\bX_i-\overline{\bX}_i)$ for $i=1,2,$ 
where $\overline{\bX}_i=[\bar{\bx}_i,...,\bar{\bx}_i]$ and $\bar{\bx}_i=\sum_{j=1}^{n_i}\bx_{j(i)}/n_i$. 
Let $\hat{\bu}_{1(i)}$ be the first (unit) eigenvector of $(\bX_i-\overline{\bX}_i)^T(\bX_i-\overline{\bX}_i)$ for $i=1,2$. 
Note that $\hat{\bu}_{1(i)}^T\bP_{n_i}=\hat{\bu}_{1(i)}^T$ when $(\bX_i-\overline{\bX}_i)^T(\bX_i-\overline{\bX}_i) \neq \bO$ for $i=1,2$. 
Then, under (A-ii) for each $\pi_i$, we have that 
\begin{equation}
\frac{\hat{\bu}_{1(1)}^T (\bX_1-\overline{\bX}_1)^T(\bX_2-\overline{\bX}_2)\hat{\bu}_{1(2)}}
{( \nu_1 \nu_2 \lambda_{1(1)}\lambda_{1(2)})^{1/2}}
=\frac{ \hat{\bu}_{1(1)}^T \sum_{s,t \le i_{\star}} \beta_{st} \bz_{os(1)}\bz_{ot(2)}^T\hat{\bu}_{1(2)} }{( \nu_1 \nu_2 \lambda_{1(1)}\lambda_{1(2)})^{1/2}}+o_p(1)
\label{A.5}
\end{equation}
as $d\to \infty$ either when $n_i$ is fixed or $n_i\to \infty$ for $i=1,2$. 
Note that $\tilde{\bh}_{1(i)}=\{\nu_i \tilde{\lambda}_{1(i)}\}^{-1/2}(\bX_i-\overline{\bX}_i)\hat{\bu}_{1(i)}$ for $i=1,2$. 
Also, note that $\bz_{os(i)}^T\bz_{os'(i)}/n_i=o_p(1)$ $(s\neq s')$ when $n_i\to \infty$ for $i=1,2$. 
Then, by combining (\ref{A.5}) with Theorem 2.1 and (\ref{A.4}), 
we can claim the result. 
\hfill $\Box$
\\[5mm]
{\it Proofs of Theorems 4.1 and 4.2}. \ 
By combining Theorem 2.1, Lemmas 2.1 and 4.1, we can claim the results. 
\hfill $\Box$
\section*{Acknowledgements}
Research of the second author was partially supported by Grant-in-Aid for Young Scientists (B), Japan Society for the Promotion of Science (JSPS), under Contract Number 26800078.
Research of the third author was partially supported by Grants-in-Aid for Scientific Research (B) and 
Challenging Exploratory Research, JSPS, under Contract Numbers 22300094 and 26540010. 




\begin{thebibliography}{99}


\bibitem[Ahn, J., Marron, J. S., Muller, K. M. and Chi, Y.-Y. (2007)]{Ahn2007}
Ahn, J., Marron, J.S., Muller, K.M., Chi, Y.-Y., 2007.
The high-dimension, low-sample-size geometric representation holds under mild conditions. 
Biometrika 94, 760-766. 

\bibitem[Aoshima, M. and Yata, K. (2011)]{Aoshima2011}
Aoshima, M., Yata, K., 2011.
Two-stage procedures for high-dimensional data. 
Sequential Anal. (Editor's special invited paper) 30, 356-399. 

\bibitem[Aoshima, M. and Yata, K. (2013)]{Aoshima2013}
Aoshima, M., Yata, K., 2013.
Asymptotic normality for inference on multisample, high-dimensional mean vectors under mild conditions. 
Methodol. Comput. Appl. Probab., in press (DOI: 10.1007/s11009-013-9370-7).

\bibitem[Armstrong et al. (2002)]{Armstrong:2002}
Armstrong, S.A., Staunton, J.E., Silverman, L.B., Pieters, R., den Boer, M.L., Minden, M.D., Sallan, S.E., 
Lander, E.S., Golub, T.R., Korsmeyer, S.J., 2002.
MLL translocations specify a distinct gene expression profile that distinguishes a unique leukemia. 
Nature Genetics 30, 41-47. 

\bibitem[Hall et al. (2005)]{Hall:2005}
Hall, P., Marron, J.S., Neeman, A., 2005.
Geometric representation of high dimension, low sample size data.
J. R. Statist. Soc. B 67, 427-444.

\bibitem[Ishii, A., Yata, K. and Aoshima, M. (2014)]{Ishii2014}
Ishii, A., Yata, K., Aoshima, M., 2014.
Asymptotic distribution of the largest eigenvalue via
geometric representations of high-dimension, low-sample-size data.
Sri Lankan J. Appl. Statist., 
Special Issue: Modern Statistical Methodologies in the Cutting Edge of Science (ed. Mukhopadhyay, N.), 81-94.

\bibitem[Jung and Marron (2009)]{Jung:2009}
Jung, S., Marron, J.S., 2009.
PCA consistency in high dimension, low sample size context. 
Ann. Statist. 37, 4104-4130. 

\bibitem[Jung et al. (2012)]{Jung:2012}
Jung, S., Sen, A., Marron, J.S., 2012. 
Boundary behavior in high dimension, low sample size asymptotics of PCA. 
J. Multivariate Anal. 109, 190-203. 

\bibitem[Srivastava, M.S. and Yanagihara, H. (2010)]{Srivastava:2010}
Srivastava, M.S., Yanagihara, H., 2010.
Testing the equality of several covariance matrices with fewer observations than the dimension. 
J. Multivariate Anal. 101, 1319-1329.



\bibitem[Yata and Aoshima (2009)]{Yata:2009}
Yata, K., Aoshima, M., 2009. 
PCA consistency for non-Gaussian data in high dimension, low sample size context. 
Commun. Statist. Theory Methods, Special Issue Honoring Zacks, S. (ed. Mukhopadhyay, N.) 38, 2634-2652.

\bibitem[Yata and Aoshima (2010)]{Yata:2010}
Yata, K., Aoshima, M., 2010.
Effective PCA for high-dimension, low-sample-size data with singular value decomposition of cross data matrix.  
J. Multivariate Anal. 101, 2060-2077.

\bibitem[Yata and Aoshima (2012)]{Yata:2012}
Yata, K., Aoshima, M., 2012.
Effective PCA for high-dimension, low-sample-size data with noise reduction via geometric representations.   
J. Multivariate Anal. 105, 193-215.

\bibitem[Yata and Aoshima (2013)]{Yata:2013}
Yata, K., Aoshima, M., 2013.
PCA consistency for the power spiked model in high-dimensional settings. 
J. Multivariate Anal. 122, 334-354.

\end{thebibliography}


\end{document}